\documentclass[12pt,a4paper]{amsart}
\newtheorem{de}{Definition}[section]

\newtheorem{re}[de]{Remark}

\newtheorem{te}[de]{Theorem}

\usepackage{soul}
\usepackage{amsfonts} 
\usepackage{amsmath}
\usepackage[all]{xy}

\newcommand{\ot}{\otg{B}}

\newcommand{\mproof}{\noindent{\bf Proof.}}
\newcommand{\twosid}[3]{\ar@<0.25ex>@{<-}[#1]^{#2}  \ar@<-1ex>[#1]_{#3}}

\def\ot{\otimes}

\begin{document}

{\bf THE UNIFICATION OF NON-ASSOCIATIVE STRUCTURES}

\bigskip

\begin{center}
 Florin F. Nichita
\end{center}

\begin{center}
 Simion Stoilow Institute of Mathematics
of the Romanian Academy\\
21 Calea Grivitei Street, 010702 Bucharest, Romania
\end{center}

%


\bigskip

\begin{abstract} 
The main non-associative algebras are Lie algebras 
and Jordan algebras.
There are several ways to unify these non-associative 
algebras and associative algebras.

\end{abstract}
\bigskip




\bigskip

\bigskip

\section{Introduction}

The main non-associative structures are Lie algebras and Jordan algebras.
Arguable less studied, Jordan algebras have applications
in physics, differential geometry, ring geometries, quantum groups,
analysis, biology, etc
(see \cite{RI}).

There are several ways to unify 
Lie algebras, Jordan algebras
 and associative algebras. 
The next section presents
structures which unify (non-)associative structures.
The last section refers to cases when the
unification of (non-)associative structures
 could be realised just in the conclusions of teorems.

This paper is related to a  communication made at the
``Workshop on Non-associative Algebras and Applications'', 
Lancaster University, UK,
July 2018.
\bigskip

All tensor products will be defined over the field $k$.

\bigskip

\section{Unification Structures}

\bigskip

\subsection{ UJLA structures} 
 The UJLA structures 
could be seen as structures which comprise the information
encapsulated in associative algebras, Lie algebras and Jordan algebras.

\begin{de}
For a $k$-space V, let 
   $ \eta : V \otimes V \rightarrow V, \ \ 
 a \otimes b \mapsto ab ,  $  be a linear map such that:
\begin{equation} \label{new}
 (ab)c + (bc)a + (ca)b  = a(bc) + b(ca) + c(ab),
\end{equation} 
\begin{equation} \label{Jordan}
 (a^2 b) a \ = \ a^2 (ba), \ \ \
 (a b) a^2 \ = \ a (b a^2), \ \ \
 (b a^2) a \ = \  (ba) a^2 , \ \ \
 a^2 (ab) \ = \  a (a^2 b),
\end{equation} 
$\forall \ a, b, c \in V $.
Then, $(V, \eta) $ is called a {\bf UJLA structure}.
\end{de}

\begin{re}
The UJLA structures  unify Jordan, Lie 
and (non-unital) associative algebras.
\end{re}
 
\begin{re} If $ (A, \ \theta )$, where $ \theta : A \otimes A \rightarrow A, \ \ 
\theta (a \otimes b) = ab $, 
 is a (non-unital) associative algebra, then we define  a UJLA structure $ (A, \ \theta' )$, where $ \theta' (a \ot b) = \alpha ab \ + \ \beta ba $, for some
 $  \ \alpha, \ \beta \in k$.
For
$\alpha = \beta = \frac{1}{2}$,  $ (A, \ \theta' )$ is a Jordan algebra, and
for $\alpha = 1 = - \beta $,  $ (A, \ \theta' )$ is a Lie algebra. 
\end{re}

\begin{te}  {\bf (Nichita \cite{upg})}
Let $(V, \eta) $ be a UJLA structure. Then,
 $(V, \eta'), \ \ \eta'(a\ot b)= [a,b]=  ab - ba $ is a Lie algebra.
\end{te}

\begin{te}  {\bf (Nichita \cite{upg})}
Let $(V, \eta) $ be a UJLA structure. Then,
 $(V, \eta'), \ \ \eta'(a\ot b)= a \circ b =  {\frac{1}{2}} (ab + ba) $ is a Jordan algebra.
\end{te}
\begin{re}
 The structures from the two above theorems are related by the relation:

$ [a, b \circ c] + [b, c \circ a] + [c, a \circ b]=0$.

\end{re}

\begin{re}
 The classification of UJLA structures is an open problem.
\end{re}

\bigskip

\subsection{Yang--Baxter equations}

The authors of \cite{a} argued that the 
Yang--Baxter equation leads to another unification of
(non-)associative structures.

\bigskip

For $ V $ a $ k$-space, we denote by
$ \   \tau : V \otimes V \rightarrow V \ot V \  $ the twist map defined by $ \tau (v \ot w) = w \ot v $, and by $ I: V \rightarrow V $
the identity map of the space $V$;
for $ \  R: V \ot V \rightarrow V \ot V  $
a $ k$-linear map, let
$ {R^{12}}= R \ot I , \  {R^{23}}= I \ot R , \
{R^{13}}=(I\ot\tau )(R\ot I)(I\ot \tau ) $.

\bigskip
\begin{de} A
{ \it  Yang-Baxter
operator} is an invertible $ k$-linear map,
$ R : V \ot V \rightarrow V \ot V $,
which satisfies the braid condition (sometimes called the {\em Yang-Baxter equation}):
\begin{equation}  \label{ybeq}
R^{12}  \circ  R^{23}  \circ  R^{12} = R^{23}  \circ  R^{12}  \circ  R^{23}.
\end{equation}
If $R$ satisfies (\ref{ybeq}) then both
$R\circ \tau  $ and $ \tau \circ R $ satisfy the {\em quantum Yang-Baxter equation} (QYBE):
\begin{equation}   \label{ybeq23}
R^{12}  \circ  R^{13}  \circ  R^{23} = R^{23}  \circ  R^{13}  \circ  R^{12}.
\end{equation}
\end{de}
Therefore, the equations (\ref{ybeq}) and (\ref{ybeq23}) are equivalent.

\bigskip

For $A$ be a (unitary) associative $k$-algebra, and $ \alpha, \beta, \gamma \in k$, the authors of
  \cite{DasNic:yan} defined the
$k$-linear map
$ R^{A}_{\alpha, \beta, \gamma}: A \ot A \rightarrow A \ot A, $
\begin{equation} \label{ybdn}
 a \ot b \mapsto \alpha ab \ot 1 + \beta 1 \ot ab -
\gamma a \ot b 
\end{equation}
which is a Yang-Baxter operator if and only if one
of the following cases holds:\\ 
(i) $ \alpha = \gamma \ne 0, \ \ \beta \ne 0 $; $ \ $
(ii) $ \beta = \gamma \ne 0, \ \ \alpha \ne 0 $; $ \ $
(iii) $ \alpha = \beta = 0, \ \ \gamma \ne 0 $.

\bigskip

  
An interesting property of (\ref{ybdn}), can be visualized in knot theory, where
the  link invariant associated to $ R^{A}_{\alpha, \beta, \gamma}$
is the Alexander polynomial.

\bigskip

For $ ( L , [,] )$  a Lie algebra over $k$,
  $ z \in Z(L) = \{ z \in L : [z,x]=0 \ \ \forall \ x \in L \} ,
  $ and $ \alpha \in k $, the authors of the papers \cite{mj} and \cite{nipo}  
defined the following Yang-Baxter operator:
${ \phi }^L_{ \alpha} \ : \ L \ot L \ \ \longrightarrow \ \  L \ot L $,
\begin{equation} \label{Lie}
x \ot y \mapsto \alpha [x,y] \ot z +  y \ot x \ .
\end{equation}

\begin{re}
 The formulas (\ref{ybdn}) and (\ref{Lie}) lead to the unification of
associative algebras and Lie algebras in the framework 
of Yang-Baxter structures. At this moment, we do not
have a satisfactory answer to the question how
 Jordan algebras fit in this framework (several partial answers were given).
\end{re}

\section{Unification of the Conclusions of Theorems}

Sometimes it is not easy to find structures which unify
theorems for  (non-)associative structures, but
we could unify just the conclusions of teorems, 
as we will see in the next theorems.

\begin{te}
 If $A$ is a Jordan algebra, a Lie algebra or an
associative algebra, and if $ \ a, b \in A$, then

$$ D: A \rightarrow A, \ \  D(x) = 
a(bx)+ b(ax)+ (ax)b- a(xb)- (xb)a- (xa)b \ \ \ \ $$

is a derivation.
\end{te}

{\bf Proof.} We consider three cases.

If $A$ is a Jordan algebra, then
 $ \ \  D(x) = 
a(bx)+ b(ax)+ (ax)b- a(xb)- (xb)a- (xa)b = a(bx)- (xa)b=
a(bx)- b(ax) $. According to \cite{b}, $D$ is a derivation.

If $A$ is a Lie algebra, then
 $ \ \  D(x) = 
a(bx)+ b(ax)+ (ax)b- a(xb)- (xb)a- (xa)b = 
a(bx)- b(ax) = a(bx)+ b(xa)= (ab)x $. So, $D$ is a derivation.

If $A$ is an associative algebra, then
 $ \ \  D(x) = 
a(bx)+ b(ax)+ (ax)b- a(xb)- (xb)a- (xa)b = (ab+ba)x - x (ab+ba)$. 
So, $D$ is a derivation.

\qed

\begin{te}
 If $A$ is a Jordan algebra, a Lie algebra or an
associative algebra, and if $ \ a, b \in A$, then
$ D: A \rightarrow A, \ \  D(x) = 
a(bx) - (xa)b \ \ \ \ $
is a derivation.
\end{te}

{\bf Proof.} We consider three cases, and follow similar steps as in the
previous proof.
\qed

\end{document}